\documentclass[12pt]{article}
\usepackage{amssymb,amsmath,amsthm}
\usepackage[utf8]{inputenc}
\usepackage[T1]{fontenc}
\usepackage{enumerate}
\usepackage[all]{xy}
\usepackage{lambda}
\def\tab{\hspace*{6mm}}
\usepackage{tgheros}
\usepackage{mathrsfs}

\usepackage{url}

\newcommand{\fin}{
\vskip 2mm
\noindent
$\Box$}

\begin{document}
\begin{center}
\LARGE Th\'{e}orie inverse de Galois sur les corps des fractions rationnelles tordus
\end{center}
\vskip 3mm
\noindent
\centerline{Angelot BEHAJAINA}
\vskip 6mm
\noindent
\begin{center}
\begin{minipage}{11cm}
{\small {\bf Abstract.---} In this article, we prove that if $H$ is a skew field of center $k$ and $\sigma$ an automorphism of finite order of $H$ such that the fixed subfield $k^{\langle \sigma \rangle}$ of $k$ under the action of $\sigma$ contains an ample field, then the inverse Galois problem has a positive answer over the skew field $H(t,\sigma)$ of twisted rational fractions.  Moreover, if $k^{\langle \sigma \rangle}$  contains either a real closed field, or an Henselian field of residue characteristic $0$ and containing all roots of unity, then the profree group of countable rank $\widehat{F}_{\omega}$ is a Galois group over $H(t,\sigma)$.\\
{\bf R\'{e}sum\'{e}.---} Dans cet article, nous montrons que si $H$ est un corps gauche de centre $k$ et $\sigma$ un automorphisme d'ordre fini de $H$ tels que le sous-corps $k^{ \langle \sigma \rangle}$ de $k$ fix\'{e} par $\sigma$ contient un corps ample, alors le probl\`{e}me inverse de Galois admet une r\'{e}ponse positive sur le corps $H(t,\sigma)$ des fractions rationnelles tordus. De plus, si  $k^{\langle \sigma \rangle}$ contient un corps qui est, soit réel clos, soit hens\'{e}lien de caract\'{e}ristique r\'{e}siduelle nulle et contenant toutes les racines de l'unité, alors le groupe prolibre de rang d\'{e}nombrable $\widehat{F}_{\omega}$ est groupe de Galois sur $H(t,\sigma)$.}
\end{minipage}
\end{center}
\vskip 5mm
\noindent

\section{Introduction}
\tab Le traditionnel problème inverse de Galois sur un corps commutatif $k$, noté {\bf $\mathrm{\bf PIG}_k$}, consiste à savoir si, {\it pour tout groupe fini $G$, il existe une extension galoisienne finie $L/k$ dont le groupe de Galois est $G$.} L'étude de ce problème dans le cas d'un corps de fractions $k(t)$ peut s'appréhender de manière géométrique si l'on demande en plus à l'extension $L/k(t)$ d'être régulière\footnote {c'est-à-dire, si $L \cap \overline{k}=k$} sur $k$. Il s'agit du problème inverse de Galois régulier, noté {\bf $\mathrm{\bf PIGR}_k$}, qui consiste lui à savoir si, {\it pour tout groupe fini $G$, il existe un $G$-revêtement $X \longrightarrow \mathbb{P}^{1}$ défini sur $k$.} Dans les années 90, Florian Pop a démontré que le {$\mathrm{ PIGR}_k$} admettait une réponse positive dès que $k$ contenait un corps ample\footnote{Un corps $k$ est dit {\it ample} si toute courbe lisse définie sur $k$ et géométriquement irréductible possède une infinité de points $k$-rationnels dès qu'elle en possède un. Par exemple, la clôture totalement réelle $\mathbb{Q}^{\mathrm{tr}}$ de $\mathbb{Q}$, le corps $\mathbb{C}((u))$ et le corps des nombres réels $\mathbb{R}$, etc. Pour les corps amples, on peut consulter aussi \cite{Jar11}, \cite{BSF13} et \cite{Pop14}.} (voir \cite{Pop96}). En fait, dans cette situation, le corps $k$ satisfait à un problème encore plus fort (voir la preuve de \cite[Proposition 10]{DL19}) :
\vskip 2 mm
\noindent
{\bf $\mathrm{\bf PIGS}_k$ (Problème Inverse de Galois Sérié sur $k$) :} {\it Pour tout groupe fini $G$, existe t-il un $G$-revêtement $X \longrightarrow \mathbb{P}^{1}$ défini sur $k$ tel que $X$ possède un point $k$-rationnel non ramifié ?}
\vskip 2mm
\noindent
De manière équivalente, ce problème consiste à savoir si, {\it pour tout groupe fini $G$, il existe une extension galoisienne finie $L/k(t)$ dont le groupe de Galois est $G$ et telle que $L$ se plonge dans le corps des séries de Laurent $k((t))$}, énoncé qui justifie la terminologie de {\it sérié} que nous utilisons pour le {$\mathrm{ PIGS}$.}
\vskip 2mm
\noindent
\tab Dans leur article \cite{DL19}, Bruno Deschamps et François Legrand s'intéressent à la problématique inverse de Galois sur des corps gauches $H$. En effet, la théorie de Galois admet une généralisation dans le cas de ces corps (voir \cite[Chapter 3]{Coh95}), qui rend légitime le {$\mathrm{PIG}_H$}. Malheureusement, la notion de clôture algébrique n'étant pas immédiate pour les corps gauches, le {$\mathrm{PIGR}_H$} n'a pas de sens. On peut néanmoins s'intéresser au {$\mathrm{PIG}_{H(t)}$} dès que l'on sait donner un analogue d'un corps de fractions sur un corps non commutatif. La théorie des anneaux de polynômes tordus de Ore en est un cadre idéal. Rappelons que, l'anneau de polynômes tordu $H[t,\sigma]$ est le $H$-espace vectoriel des polynômes sur $H$ muni du produit qui vérifie $ta=\sigma(a)t$ pour tout $a \in H$. Cet anneau possède un unique corps de fractions $H(t,\sigma)$ (voir \cite[Chapter 2]{Coh95}). L'un des résultats principaux de [DL19] affirme que, si $H$ est de dimension finie sur son centre $k$ et si $\sigma=\mathrm{Id}$, alors 
$$  \mathrm{PIGS}_k \Longrightarrow \mathrm{PIG}_{H(t,\mathrm{Id})} . $$ 
\vskip 2mm
\noindent 
Leur méthode consiste à contrôler, grâce à l'étude de la norme réduite, l'extension des scalaires à $H$. De ce point de vue, la condition de finitude de $H$ sur $k$ est essentielle. Une application de ce résultat a été récemment donnée dans \cite{ALP20}, où il est montré que l’anneau $H(X)$ des fonctions polynomiales en la variable X et à coefficients dans $H$ est isomorphe à un certain corps $H'(t,\mathrm{Id})$. Les auteurs en déduisent, par \cite{DL19}, que le problème inverse de Galois admet une réponse sur $H(X)$ dès que $k$ contient un corps ample. 
\vskip 2mm
\noindent
\tab Dans cet article, tout en continuant à utiliser la méthode de l'extension des scalaires, nous introduisons une nouvelle approche utilisant les corps de séries de Laurent tordus, qui nous permet de montrer la généralisation du résultat de Deschamps--Legrand suivante :
\vskip 2mm
\noindent
{\bf Th\'{e}or\`{e}me A.---} {\it Soient $H$ un corps gauche de centre $k$ ($H$ non nécessairement de dimension finie sur $k$) et $\sigma$ un automorphisme d'ordre fini de $H$. Si $k^{\langle \sigma \rangle}$ désigne le corps des invariants de $k$ par $\sigma$, alors on a
$$  \mathrm{PIGS}_{k^{\langle \sigma \rangle}} \Longrightarrow \mathrm{PIG}_{H(t,\sigma)}. $$ 
En particulier, si $k^{\langle \sigma \rangle}$ contient un corps ample, alors le $\mathrm{PIG}_{H(t,\sigma)}$ admet une r\'{e}ponse positive.}
\vskip 2mm
\noindent
Un exemple d'application du théorème A ne rentrant pas dans le cadre des travaux de \cite{DL19} peut être donné en considérant le corps $H=k(u,\sigma)$ où $k/k_{0}$ désigne une extension galoisienne vérifiant $\langle \sigma \rangle=\mathrm{Gal}(k/k_{0})$ infini et $k_{0}$ ample (e.g. $k_{0}=\mathbb{C}((x))$ et $k=\overline{\mathbb{C}((x))}$). Voir les exemples d'application dans le \S 2.1 pour plus de détails.
\vskip 2mm
\noindent
\tab La théorie de Galois admet aussi une généralisation au cas infini. Dans \cite[ Chapter VII : \S6 ]{Jac56}, il est montré que si $\Omega/H$ est une extension galoisienne extérieure et \textit{relativement algébrique à gauche} \footnote {dans le sens où pour tout élément $x \in \Omega$, le sous-corps $H(x)$ de $M$ engendré par $H$ et $x$ est de dimension finie à gauche sur $H$. Cette terminologie apparaît dans \cite[\S6 : Definition 1]{Jac56}.}, le groupe de Galois $\mathrm{Gal}(\Omega/H)$ de $\Omega/H$ est profini. Plus précisément, ce groupe est la limite projective des groupes de Galois $\mathrm{Gal}(M/H)$ des sous-extensions galoisiennes finies $M/H$ de $\Omega/H$. Par ailleurs, la théorie de Krull reste valable et il y a la correspondance entre les corps intermédiaires et les sous-groupes fermés (voir \cite[\S6 Theorem 1]{Jac56}). Dans cette perspective, on s'inspire des travaux de Pierre D\`{e}bes et Bruno Des\-champs sur la $\psi$-liberté des corps (présentés dans \cite{DD04}), pour nous intéresser \`{a} la possibilit\'{e} de r\'{e}aliser le groupe prolibre $\widehat{F}_{\omega}$ comme groupe de Galois d'une extension de $H(t,\sigma,\delta)$. A cet effet, nous montrons le
\vskip 2mm
\noindent
{\bf Th\'{e}or\`{e}me B.---} {\it Soient $H$ un corps gauche de centre $k$ et $\sigma$ un automorphisme d'ordre fini de $H$. Supposons que le corps des invariants $k^{\langle \sigma \rangle}$ contient un corps $k_{0}$ qui est, soit réel clos, soit hens\'{e}lien de caract\'{e}ristique r\'{e}siduelle nulle et contenant toutes les racines de l'unité. Alors, il existe une extension $L/H(t,\sigma)$ galoisienne extérieure, relativement algébrique à gauche et de groupe de Galois $\widehat{F}_{\omega}$.}
\vskip 2mm
\noindent
Ce théorème montre, par exemple, que le groupe $\widehat{F}_{\omega}$ est groupe de Galois sur $\mathbb{H}(t)$, le corps des fractions rationnelles tordu à indéterminée centrale, où $\mathbb{H}$ est le corps des quaternions de Hamilton (voir l'exemple d'application dans le \S 2.2).
\vskip 2mm
\noindent
\vskip 2mm
\noindent
{\bf Remarque.---} Dans certaines situations, les corps de fractions rationnelles tordus avec dérivation $H(t,\sigma,\delta)$ peuvent se ramener au cas où $\delta=0$, ce qui permet d'appliquer les résultats de ce travail. En effet, il est par exemple montré dans \cite[Theorem 2.3.1]{Coh95} que  si $\delta \neq 0$ et $k^{\langle \sigma \rangle} \neq k$, alors les corps $H(t,\sigma,\delta)$ et $H(t,\sigma)$ sont isomorphes. Nous n'avons pas inclus ces cas dans les théorèmes A et B pour ne pas en alourdir artificiellement les énoncés.
\vskip 2mm
\noindent
{\bf Remerciements.---} Je tiens à remercier mes directeurs de thèse, Bruno Deschamps et François Legrand, pour leurs nombreuses relectures, commentaires utiles et précieuses suggestions. Je remercie également le GDRI GANDA pour son soutien financier.

\section{Preuves des th\'{e}or\`{e}mes }\label{sec3}
\subsection{Théorème A}\label{thA}
\tab Commençons par rappeler la stratégie de l'extension galoisienne des scalaires présentée dans \cite{DL19} : si l'on se donne un corps gauche $K$ et un corps $k$ contenu dans le centre $Z(K)$ de $K$, alors on a le 
\vskip 2mm
\noindent
{\bf Lemme 2.1.1.---} Si $L/k$ désigne une extension algébrique galoisienne de groupe $G$ telle que la $k$-algèbre $M=K \otimes_{k} L$ soit un corps, alors :
\vskip 1mm
\noindent
1) Si l'on considère $K$ comme sous-corps de $M$, l'application 
\begin{equation}\label{phidef}
\begin{array}{ccccccccc}
\Psi & : & \mathrm{Gal}(L/k) & \rightarrow & \mathrm{Aut}(M/K)     &  &           &                    &  \\
             &   &   \sigma         & \mapsto     &  \Psi(\sigma) &: & M & \rightarrow & M \\
             &   &                      &                   &                             &  & h\otimes x & \mapsto & h \otimes \sigma(x) 
\end{array}
\end{equation} est un isomorphisme de groupes.
\vskip 1mm
\noindent
2) L'extension $M/K$ est galoisienne de groupe $G$.
\vskip 2mm
\noindent
{\bf Preuve :} 1) On voit facilement que $\Psi$ est un morphisme de groupes. Il reste à montrer qu'elle est bijective. L'injectivité est immédiate. Pour la surjectivité, considérons $\eta \in \mathrm{Gal}(M/K)$. On remarque tout d'abord que $Z(M)$ est corps. Vu que $\eta(Z(M))=Z(M)$, la restriction de $\eta$ \`{a} $Z(M)$ est un $k$-automorphisme. Or $(k\otimes_{k} L)/k$ est une sous-extension galoisienne de $Z(M)/k$, donc $\eta(k \otimes_{k}L)=k \otimes_{k}L$. Ainsi, on peut définir une application $\sigma_{\eta}:L \rightarrow L$ qui à $l \in L$ associe l'unique $\tilde{l} \in L$ v\'{e}rifiant $\eta(1 \otimes l)=1 \otimes \tilde{l}$. On vérifie aisément que $\sigma_{\eta}$ est un isomorphisme d'anneaux, donc $\sigma_{\eta} \in \mathrm{Gal}(L/k)$. Il reste \`{a} v\'{e}rifier que l'on a $\Psi(\sigma_{\eta})=\eta$. A cet effet, pour $h \in K$ et $l \in L$, on a $\Psi(\sigma_{\eta})(h \otimes l)=\Psi(\sigma_{\eta})(h \otimes 1)\Psi(\sigma_{\eta})(1 \otimes l)=(h\otimes 1) (1 \otimes \sigma_{\eta}(l))=\eta(h \otimes 1) \eta (1 \otimes l)= \eta ((h \otimes 1)(1 \otimes l))=\eta (h \otimes l) 
$. Ainsi, on a $\Psi(\sigma_{\eta})=\eta$.
\vskip 1mm
\noindent
2) On voit immédiatement que la sous-algèbre des invariants de $M$ par l'action de $\mathrm{Aut}(M/K)$ vaut $K$, ce qui conclut. 
\fin
\vskip 2mm
\noindent
\tab 
La difficulté dans la méthode de l'extension galoisienne des scalaires consiste à assurer que l'algèbre $M=K \otimes_k L$ reste un corps. Dans \cite{DL19}, les auteurs y arrivent en considérant pour $K$ un corps de fractions tordus à indéterminée centrale dont le corps des constantes $H$ est de dimension finie sur $k$. Ils contrôlent alors les zéros de la norme réduite dans l'extension $L/k$. Ceci limite de manière rédhibitoire l'étude au cas où $K=H(t,\mathrm{Id})$ avec $[H:k]<+ \infty$. La nouvelle approche que nous proposons pour le cas $K=H(t,\sigma)$ consiste à plonger $K$, $L$ et $M$ dans un même corps (de séries). Ceci va nous permettre d'étendre le résultat de \cite{DL19} à des cas où $\sigma \neq \mathrm{Id}$ et $[H:k]=+\infty$.
\vskip 2mm
\noindent
\tab Etant donnés un corps gauche $H$ de centre $k$ et $\sigma$ un automorphisme de $H$, rappelons que le corps des séries de Laurent tordu, noté $H((t,\sigma))$, est constitué des séries formelles $\sum_{i \geq i_{0}}x_{i}t^{i}$ ($i_{0} \in \mathbb{Z}$ et $x_{i} \in H$) et est muni du produit $ta=\sigma(a)t$ pour tout $a \in H$. Le corps $H(t,\sigma)$ se plonge canoniquement dans $H((t,\sigma))$ par $t \mapsto t$ et $a \mapsto a$ ($a \in H$). Si l'on suppose que $\sigma$ est d'ordre $n \geq 1$, alors le corps $k^{\langle \sigma \rangle}(t^{n})$ est à la fois dans le centre de $H(t,\sigma)$ et dans celui de $H((t,\sigma))$. Pour plus de détails, nous renvoyons le lecteur à \cite[Chapter 2]{Coh95}. 
\vskip 2mm
\noindent
{\bf Proposition 2.1.2.---} {\it Soient $H$ un corps gauche de centre $k$ (de dimension non nécessairement finie sur $k$) et $\sigma$ un automorphisme de $H$ d'ordre $n \geq 1$. Si $k^{\langle \sigma \rangle}$ désigne le corps des invariants de $k$ par l'action de $\sigma$, alors l'application $k^{\langle \sigma \rangle}(t^{n})$-linéaire
$$\tau: H(t,\sigma)\otimes_{k^{\langle \sigma \rangle}(t^{n})}k^{\langle \sigma \rangle}((t^{n}))  \rightarrow H((t,\sigma)),$$ qui envoie chaque tenseur $x \otimes l$ sur $xl$, est un monomorphisme de $k^{\langle \sigma \rangle}(t^{n})$-alg\`{e}bres.
\vskip 2mm
\noindent
\tab En particulier, pour toute extension intermédiaire $k^{\langle \sigma \rangle}(t^{n}) \subset L \subset k^{\langle \sigma \rangle}((t^{n}))$ de dimension finie sur $k^{\langle \sigma \rangle}(t^{n})$, l'algèbre $H(t,\sigma)\otimes_{k^{\langle \sigma \rangle}(t^{n})}L$ est un corps.}
\vskip 2mm
\noindent
{\bf Preuve :} En toute généralité, si $A, B, C$ désignent trois $k$-algèbres et $f:A \rightarrow C$, $g: B \rightarrow C$ sont deux morphismes de $k$-algèbres, alors l'application $k$-linéaire de $A\otimes_{k}B \longrightarrow C$ qui au tenseur $a \otimes b$ associe $f(a)g(b)$, est un morphisme d'algèbres si et seulement si $\mathrm{Im}(g)$ est inclus dans le {\it commutant} \footnote{Etant donné un anneau $R$ et une partie $S$ de $A$, le commutant de $S$ dans $R$ est l'ensemble des $x \in R$ qui commutent avec tous les éléments de $S$.} dans $C$ de $\mathrm{Im}(f)$ (ce qui équivaut à dire que $\mathrm{Im}(f)$ est inclus dans le commutant dans $C$ de $\mathrm{Im}(g)$). Dans notre situation, $k^{\langle \sigma \rangle}((t^{n}))$ est inclus dans le centre de $H((T,\sigma))$, ce qui montre que $\tau$ est bien un morphisme d'algèbres.
\vskip 2mm
\noindent
\tab Pour l'injectivité de $\tau$, considérons $h_{1},\cdots,h_s \in  H(t,\sigma)$ et $z_{1},\cdots,z_{s} \in k^{\langle \sigma \rangle}((t^{n}))$ tels que $x=\sum_{m=1}^{s} h_{m} \otimes z_{m}$ vérifie $\tau(x)=0$. Puisque $H[t,\sigma]$ est un anneau de Ore (voir \cite[Chapter 2]{Coh95}), en appliquant $s$ fois la possible écriture $b^{-1}a=a'b'^{-1}$ ($a,b,a',b' \in H[t,\sigma]$), on déduit l'existence d'un polynôme $u \in H[t,\sigma]-\{0\}$ tel que pour tout $m \in \{1, \cdots, s \}$, on a $h_{m}u \in H[t,\sigma]$. Comme $u \otimes 1$ est inversible et comme $\tau(x)\tau(u\otimes 1)=\tau(\sum_{m=1}^{s} h_{m}u \otimes z_{m})$, pour montrer que $x=0$, on peut supposer que les $h_m$ sont dans $H[t,\sigma]$. Le même type de manipulation permet de supposer que les $z_m$ sont dans $k^{\langle \sigma \rangle}[[t^n]]$. 
\vskip 2mm
\noindent
\tab Fixons une $k^{\langle \sigma \rangle}$-base $(e_{j})_{ j \in J}$ de $H$.  En utilisant la centralité de $k^{\langle \sigma \rangle}(t^n)$, on peut écrire chacun des $h_m$ sous la forme d'une somme finie 
$$\begin{array}{lll}
\displaystyle h_m(t)& = & \displaystyle \sum_{l=0}^{l_m} t^{ln}\sum_{k=0}^{n-1} a_{l,k,m}t^{k}\ \hbox{\rm\small avec}\ a_{l,k,m}\in H \\
&=& \displaystyle\sum_{l=0}^{l_m} t^{ln}\sum_{k=0}^{n-1} \left(\sum_{j\in J_m}\lambda_{l,k,j,m}e_j\right) t^{k}\ \hbox{\rm\small avec}\ J_m\subset J\ \hbox{\rm\small un sous-ensemble fini choisi}\\
&& \hbox{\rm\small uniformément pour}\  l\ \hbox{\rm\small et}\ k\ \hbox{\rm\small variant et}\ \lambda_{l,k,j,m}\in k^{\langle \sigma \rangle} \\
&=& \displaystyle \sum_{k=0}^{n-1}\sum_{j\in J_m}e_jt^kP_{k,j,m}(t^n)\ \hbox{\rm\small avec}\ P_{k,j,m}(t^n)\in k^{\langle \sigma \rangle}[t^n].\\
\end{array}$$
On a alors
$$h_m\otimes z_m=\sum_{k=0}^{n-1}\sum_{j\in J_m}(e_jt^kP_{k,j,m}(t^n))\otimes z_m=\sum_{k=0}^{n-1}\sum_{j\in J_m}(e_jt^k)\otimes (P_{k,j,m}(t^n)z_m)$$
et en sommant pour $m=1,\cdots ,s$, on trouve finalement une partie $J_0=\bigcup_m J_m\subset J$ finie et des séries 
$$z_{k,j}=\sum_{m=1}^s P_{k,j,m}(t^n)z_m\in  k^{\langle \sigma \rangle}[[t^n]]$$
(en ayant pris soin de poser $P_{k,j,m}(t^n)=0$ si $j\notin J_m$) telles que
$$x=\sum_{m=1}^s h_m\otimes z_m=\sum_{k=0}^{n-1}\sum_{j\in J_0}(e_jt^k)\otimes z_{k,j}.$$
En posant formellement $\displaystyle z_{k,j}=\sum_{l\geq 0}\alpha_{k,j,l}t^{ln}$ (avec $\alpha_{k,j,l} \in k^{\langle \sigma \rangle}$) et en utilisant la centralité, on trouve alors
$$\tau(x)=\sum_{k=0}^{n-1} \sum_{j\in J_0}e_jt^k\sum_{l\geq 0}\alpha_{k,j,l}t^{ln}=\sum_{l\geq 0}\sum_{k=0}^{n-1}\left(\sum_{j\in J_0}\alpha_{k,j,l}e_j\right)t^{ln+k}=0.$$
On en déduit que, pour tout $k=0,\cdots, n-1$ et tout $l\geq 0$, on a $\sum_{j\in J_0}\alpha_{k,j,l}e_j=0$. La famille $(e_j)_j$ étant une $k^{\langle \sigma \rangle}$-base de $H$, on a alors $\alpha_{k,j,l}=0$ pour tous $k,j,l$. Ainsi, toutes les séries $z_{k,j}$ sont nulles et donc $x=0$.
\vskip 2mm
\noindent
\tab Le dernier point de la proposition 2.1.2. découle alors du 
\vskip 2mm
\noindent
{\bf Lemme 2.1.3.---} \textit{Si $A$ désigne un anneau unitaire intègre contenant un corps $K$ (non nécessairement commutatif) tel que $A$ soit de dimension finie en tant que $K$-espace vectoriel à gauche (pour la multiplication à droite dans $A$), alors $A$ est un corps.}
\vskip 2mm
\noindent
{\bf Preuve du lemme 2.1.3. :}  Fixons $x \in A \setminus \{0\}$ et notons $m_x : a \mapsto ax$ la multiplication à droite par $x$. Il est clair que $m_x \in \mathscr{L}_K(A)$, mais comme $A$ est supposé intègre, cette application linéaire est injective. Rappelons que les théorèmes classiques (théorème du rang, etc.) d'algèbre linéaire sur les corps commutatifs restent valables sur les corps gauches (voir \cite[\S 3.11]{Ehr11} et \cite[\S 3.12]{Ehr11} ). En application du théorème du rang et puisque $A$ est supposé de dimension finie sur $K$, on en déduit que $m_x$ est surjective. Ainsi, il existe $x^{-1}\in A$ tel que $m_x(x^{-1})=x^{-1}x=1$. Ainsi, tout $x \neq 0$ possède un inverse à gauche mais comme $1$ est un neutre bilatère, finalement, tout $x \neq 0$ possède un inverse à droite.
\fin
\vskip 2mm
\noindent
En effet, considérons $k^{\langle \sigma \rangle}(t^{n}) \subset L \subset k^{\langle \sigma \rangle}((t^{n}))$ tel que le corps $L$ soit de dimension finie sur $k^{\langle \sigma \rangle}(t^{n})$. Puisque $H(t,\sigma) \otimes_{k^{\langle \sigma \rangle}(t^{n})} L$ se plonge dans le corps $H((t,\sigma))$, c'est un anneau unitaire intègre. En tant que $H(t,\sigma)$-espace vectoriel à gauche, $H(t,\sigma) \otimes_{k^{\langle \sigma \rangle}(t^{n})} L$ est de dimension $[L:k^{\langle \sigma \rangle}(t^{n})] < +\infty$. Le lemme 2.1.3. montre alors que $H(t,\sigma) \otimes_{k^{\langle \sigma \rangle}(t^{n})} L$ est bien un corps.
\fin
\vskip 2mm
\noindent
\tab Maintenant, venons-en à la
\vskip 2mm
\noindent
{\bf Preuve du théorème A :} Supposons que le $\mathrm{PIGS}_{k^{\langle \sigma \rangle}}$ admet une réponse positive. Soit $G$ un groupe fini. En notant $n \geq 1$ l'ordre de $\sigma$, par hypothèse, il existe une extension galoisienne finie $L/k^{\langle \sigma \rangle}(t^{n})$ de groupe $G$ telle que $L \subset k^{\langle \sigma \rangle}((t^{n}))$. Par la proposition 2.1.2, l'algèbre $H(t,\sigma) \otimes_{k^{\langle \sigma \rangle}(t^{n})}L$ reste un corps. De plus, par le lemme 2.1.1, l'extension $H(t,\sigma) \otimes_{k^{\langle \sigma \rangle}(t^{n})}L/H(t,\sigma)$ est galoisienne de groupe de Galois $G$. Par conséquent, le $\mathrm{PIG}_{H(t,\sigma)}$ admet une réponse positive. En particulier, si $k^{\langle
\sigma \rangle}$ contient un corps ample, alors, comme rappelé dans l'introduction, le $\mathrm{PIGS}_{k^{\langle \sigma \rangle}}$ admet une réponse positive et il en est donc de même du $\mathrm{PIG}_{H(t, \sigma)}$.
\fin
\vskip 2mm
\noindent
{\bf Remarque.---}{Soient $H$ un corps gauche de centre $k$ et $\sigma$ un automorphisme d'ordre fini $n \geq 1$ de $H$. Soit $G$ un groupe fini vérifiant l'un des quatre cas suivants :
\vskip 1mm
\noindent
1) $G$ abélien et $k$ infini,
\vskip 0.5mm
\noindent
2) $G=S_n$ ($n \geq 3$) et $k$ infini,
\vskip 0.5mm
\noindent
3) $G=A_n$ ($n \geq 4$) et $k$ est de caractéristique nulle,
\vskip 0.5mm
\noindent
4) $G$ résoluble et $k$ de caractéristique positive.
\vskip 1mm
\noindent
Par la preuve de \cite[Proposition 12]{DL19}, il existe $k_0 \subset k^{\langle \sigma \rangle}$ et une extension galoisienne $L/k_0(t^n)$ de groupe $G$ et vérifie $L \subset k_0((t^{n}))$. En conséquence, $Lk^{\langle \sigma \rangle}/k^{\langle \sigma \rangle}(t^n)$ est galoisienne de groupe de $G$ vérifiant $Lk^{\langle \sigma \rangle} \subset k^{\langle \sigma \rangle}((t^n))$. Avec la même preuve du théorème A, $G$ est groupe de Galois sur $H(t,\sigma)$. }
\vskip 2mm
\noindent
{\bf Exemples d'application.---} Dans \cite{DL19}, les auteurs fournissent des exemples de corps non commutatifs, de dimensions finies sur leurs centres et sur lesquels le {$\mathrm{PIG}$} admet une réponse positive. Ils montrent que si $H$ est un corps gauche contenant un corps ample dans son centre, alors le {$\mathrm{PIG}_{H(t)}$} admet une réponse positive, où $H(t)$ est de dimension finie sur son centre.  Notre théorème A est bien plus général comme le montre les exemples suivants.
\vskip 0.5mm
\noindent
a) En considérant la conjugaison complexe $\tau$ sur $\mathbb{C}$, le corps $\mathbb{C}(t,\tau)$ est bien une algèbre simple centrale de centre $\mathbb{R}(t^2)$, qui contient un corps ample. Le corps $\mathbb{C}(t, \tau)$ est bien de dimension finie sur son centre mais les résultats de [DL19] ne peuvent s'appliquer car l'indéterminée n'est pas centrale. Notre théorème A assure cependant une réponse positive au {$\mathrm{PIG}_{\mathbb{C}(t,\tau)}$}.   
\vskip 1mm
\noindent
b) Lorsque $k_0$ est un corps ample à groupe de Brauer nul (par exemple lorsque $k_0$ est PAC\footnote{Un corps $k$ est dit {\it PAC} (Pseudo Algébriquement Clos) si toute variété non vide définie sur $k$ admet au moins un $k$-point rationnel. Par exemple, un corps algébriquement clos. Le fait qu'être PAC implique avoir un groupe de Brauer nul vient de \cite[Theorem 11.6.4]{FJ08}. Voir \cite{FJ08} pour plus de détails sur les corps PAC.}), les résultats de \cite{DL19} sont sans intérêt puisqu'il n'existe pas d'algèbre simple centrale non triviale de centre $k_0$. Pour autant, notre théorème A s'applique non trivialement à des corps $H$ de centre $k_0$ un corps ample à groupe de Brauer nul. Par exemple, si l'on prend $k_0 =\mathbb{C}((u))$, on a $\overline{k_0}=\mathrm{Puis}(\mathbb{C})$ et  $\mathrm{Gal}(\overline{k_0}/k_0) \simeq \widehat{\mathbb{Z}}$. En considérant $\sigma$ un générateur de $\widehat{\mathbb{Z}}$, on voit que $H=\overline{k_{0}}(y,\sigma)$ est de centre $k_0$ et le théorème A assure que le $\mathrm{PIG}$ est vrai sur $H(t)$. 
\subsection{Théorème B}
\tab Dans leur article \cite{DD04}, les auteurs introduisent la notion galoisienne de $\psi$-liberté, dont la définition peut s'adapter au cas des corps gauches. Rappelons en quelques mots de quoi il s'agit. Un corps commutatif $H$ est dit \textit{$\psi$-libre} si pour tout syst\`{e}me projectif $(G_{n},s_{n})_{n \in \mathbb{N}}$ de groupes finis où $s_{n}:G_{n} \twoheadrightarrow G_{n-1}$ est un épimorphisme pour tout $n \geq 1$ (ce que Dèbes et Deschamps appellent dans leur article un système complet de groupes finis), il existe une tour d'extensions finies galoisiennes $(H_{n}/H)_{n \in \mathbb{N}}$ et une suite d'isomorphismes $(\epsilon_{n} : G_{n} \rightarrow \mathrm{Gal}(H_{n}/H) )_{n \in \mathbb{N}}$ tels que pour tout $n\in \mathbb{N}$, on a  
\[
\epsilon_{n} \circ s_{n+1}  = \mathrm{res}^{H_{n+1}/H}_{H_{n}/H} \circ \epsilon_{n+1}
,\] où $\mathrm{res}^{H_{n+1}/H}_{H_{n}/H}: \mathrm{Gal}(H_{n+1}/H) \rightarrow \mathrm{Gal}(H_{n}/H)$ désigne l'épimorphisme de restriction.
\vskip 2mm
\noindent 
Dans ce cas, on dit que la tour $(H_{n}/H)_{n \in \mathbb{N}}$ {\it r\'{e}alise} le syst\`{e}me complet $(G_{n},s_{n})_{n \in \mathbb{N}}$. Cette définition n'est pas directement applicable dans le cas où $H$ est un corps gauche. En effet (voir \cite[Chapter 3]{Coh95}), si l'on dispose d'une tour d'extensions galoisiennes finies $M /L / H$, la traditionnelle suite exacte
\[
1 \rightarrow \mathrm{Gal}(M/L) \rightarrow \mathrm{Gal}(M/H) \rightarrow \mathrm{Gal}(L/H) \rightarrow 1
,\] valable dans le cas commutatif, doit être remplacée par la suite
\[
1 \rightarrow \mathrm{Gal}(M/L) \rightarrow \mathscr{N}_{\mathrm{Gal}(M/H)}(\mathrm{Gal}(M/L)) \rightarrow \mathrm{Gal}(L/H) \rightarrow 1,
\] où $\mathscr{N}_{\mathrm{Gal}(M/H)}(\mathrm{Gal}(M/L))$ désigne le normalisateur de $\mathrm{Gal}(M/L)$ dans $\mathrm{Gal}(M/H)$ (dans cette situation, on sait alors que $\mathrm{Gal}(M/H)$ est le plus petit $N$-sous-groupe \footnote{Un sous-groupe $J$ de $\mathrm{Gal}(M/H)$ est dit {\it $N$-groupe } si l'ensemble $\{x \in M \setminus \{0 \} \mid \mathrm{I}(x) \in J \} \cup \{0 \}$, où $I(x)(m)=x m x^{-1}$ pour tout $m \in M$, est une algèbre sur le centre $Z(M)$ de $M$.} contenant le normalisateur). Ainsi, l'application $\mathrm{res}^{M/H}_{L/H}$ n'est pas nécessairement définie. Une manière de pallier à ce problème consiste à considérer des extensions {\it extérieures}. En effet, dans cette situation, les éléments $x \in M$ tels que l'automorphisme intérieur $I(x): m \in M \mapsto x m x^{-1} \in M$ appartiennent au groupe $\mathscr{N}_{\mathrm{Gal}(M/H)}(\mathrm{Gal}(M/L))$ sont exactement les éléments du centre $Z(M)$ de $M$. Puisque $Z(M)$ est un corps, le groupe $\mathscr{N}_{\mathrm{Gal}(M/H)}(\mathrm{Gal}(M/L))$ est donc un $N$-groupe et l'on a donc $\mathscr{N}_{\mathrm{Gal}(M/H)}(\mathrm{Gal}(M/L))=\mathrm{Gal}(M/H)$. On peut donc définir la $\psi$-liberté d'un corps gauche de la même manière que dans le cas commutatif, mais en demandant que chaque extension $H_n/H$ soit extérieure. Un autre intérêt à ne considérer que des extensions extérieures est que, dans cette situation, la théorie de Krull des extensions galoisiennes infinies reste valide dans le cas relativement algébrique, comme mentionné dans l'introduction. Pour ces raisons, nous étendons la définition de $\psi$-liberté comme suit :
\vskip 2mm
\noindent
{\bf Définition 2.2.1.---}{\it Un corps $H$ est dit {\it $\psi$-libre} si pour tout système complet de groupes finis $(G_n,s_n)_{n \in \mathbb{N}}$, il existe une tour d'extensions $(H_n/H)_{n \in \mathbb{N}}$ telle que :
\vskip 1mm
\noindent
a) $\Omega =\bigcup_{n} H_n$ est une extension galoisienne extérieure, relativement algébrique à gauche,
\vskip 0.5 mm
\noindent
b) pour tout $n \in \mathbb{N}$, $H_n/H$ est galoisienne (nécessairement extérieure) de groupe $G_n$ et il existe $\epsilon_{n}:G_{n} \simeq \mathrm{Gal}(H_n/H)$,
\vskip 0.5mm
\noindent vérifiant  
\[
\epsilon_{n} \circ s_{n+1}  = \mathrm{res}^{H_{n+1}/H}_{H_{n}/H} \circ \epsilon_{n+1}
\] pour tout $n \in \mathbb{N}$, où $\mathrm{res}^{H_{n+1}/H}_{H_{n}/H}: \mathrm{Gal}(H_{n+1}/H) \rightarrow \mathrm{Gal}(H_{n}/H)$ désigne l'épimorphisme de restriction.}
\vskip 2mm
\noindent
\tab Eu égard aux propriétés de la théorie de Galois infinie des corps gauches, la preuve présentée dans \cite[Proposition 1.2]{DD04} montre immédiatement que $\psi$-liberté d'un corps gauche équivaut juste à l'existence d'une extension relativement algébrique galoisienne extérieure de groupe de Galois le groupe prolibre de rang dénombrable $\widehat{F}_{\omega}$ (c'est-à-dire à la réalisation d'un système complet de groupes finis particulier qui définit ce groupe).
\vskip 2mm
\noindent
{\bf Proposition 2.2.2.---} {\it
Soient $H$ un corps gauche, $k$ un corps contenu dans le centre de $H$ et $(L_{n}/k)_{n \in \mathbb{N}}$ une tour d'extensions galoisiennes qui r\'{e}alise un syst\`{e}me complet de groupes finis $(G_{n},s_{n})_{n \in \mathbb{N}}$.  Supposons que $M_{n}=H \otimes_{k} L_{n}$ reste un corps pour tout $n\in \mathbb{N}$. Alors :
\begin{enumerate}[1)]
\item pour tout $n \in \mathbb{N}$, l'extension $M_{n}/H$ est galoisienne extérieure relativement algébrique à gauche de groupe $G_{n}$,  
\item si l'on pose $L=\bigcup_{n \in \mathbb{N}}L_{n}$, alors $M=L \otimes_{k}H=\bigcup_{n \in \mathbb{N}}M_{n}$ est un corps qui est une extension galoisienne extérieure relativement algébrique à gauche de $H$, 
\item la tour d'extensions galoisiennes $(M_{n}/H)_{n \in \mathbb{N}}$ r\'{e}alise le syst\`{e}me $(G_{n},s_{n})_{n \in \mathbb{N}}$ de sorte que 
\[
\mathrm{Gal}(M/H) \simeq \lim\limits_{\substack{\longleftarrow \\ n}}G_{n}.
\] 
\end{enumerate}}
\vskip 2mm
\noindent
{\bf Preuve :} 1) Soit $n \in \mathbb{N}$. Comme $M_{n}$ est un corps, l'extension $M_{n}/H$ est galoisienne de groupe de Galois $G_{n}$. En remarquant $[M_{n}:H]=\mathrm{Gal}(M_{n}/H)$, on déduit de \cite[Théorème]{Des18} que $M_{n}/H$ est extérieure. Le fait que $M_n/H$ soit relativement algébrique est clair.
\vskip 2mm
\noindent
2) La suite $(M_{n})_{n \in \mathbb{N}}$ filtre $M$, donc $M$ est un corps. D'une part, l'extension $M/H$ est galoisienne car $L/k$ l'est. D'autre part, $M/H$ est extérieure car tout $x \in M$ est dans un certain $M_{n}$. De plus, l'extension $M/H$ est relativement algébrique en étant filtrée par les extensions relativement algébriques $M_n/H$ ($n \in \mathbb{N}$).
\vskip 2mm
\noindent
3) Par hypothèse, il existe une suite d'isomorphismes $(\epsilon_{n} : G_{n} \rightarrow \mathrm{Gal}(L_{n}/k))_{n \in \mathbb{N}}$ telle que pour tout $n\in \mathbb{N}$, on a
\begin{equation}
\epsilon_{n} \circ s_{n+1}= \mathrm{res}^{L_{n+1}/k}_{L_{n}/k} \circ \epsilon_{n+1}, \label{ddd}
\end{equation} o\`{u} $\mathrm{res}^{L_{n+1}/k}_{L_{n}/k}$ d\'{e}signe le morphisme de restriction. Pour tout $n \in \mathbb{N} $, en considérant l'isomorphisme $\Psi_{n}: \mathrm{Gal}(L_{n}/k) \rightarrow \mathrm{Gal}(M_n/H)$ défini dans \eqref{phidef}, on vérifie que l'on a
\begin{equation}
\Psi_{n} \circ \mathrm{res}^{L_{n+1}/k}_{L_{n}/k}=\mathrm{res}^{M_{n+1}/H}_{M_{n}/H} \circ \Psi_{n+1} \label{equaas}
\end{equation} pour tout $n \in \mathbb{N}$. Pour tout $n \in \mathbb{N}$, on pose $\tilde{\epsilon}_{n}=\Psi_{n} \circ \epsilon_{n}$. Alors, on a
\begin{align*}
 \tilde{\epsilon}_{n} \circ s_{n+1} =  \Psi_{n} \circ \left( \epsilon_{n}  \circ s_{n+1} \right)  = \left(\Psi_{n} \circ \mathrm{res}^{L_{n+1}/k}_{L_{n}/k} \right) \circ \epsilon_{n+1}&= \mathrm{res}^{M_{n+1}/H}_{M_{n}/H} \circ \left(\Psi_{n+1} \circ \epsilon_{n+1} \right) \\
&= \mathrm{res}^{M_{n+1}/H}_{M_{n}/H} \circ \tilde{\epsilon}_{n+1}
 \end{align*}pour tout $n \in \mathbb{N}$. La deuxi\`{e}me (resp. troisi\`{e}me) \'{e}galit\'{e} vient de l'\'{e}galit\'{e} \eqref{ddd} (resp. \eqref{equaas}). Il reste à montrer l'isomophisme. Pour tout $n \in \mathbb{N}$, on a $\mathrm{res}^{M/H}_{M_{n+1}/H} \circ \mathrm{res}^{M_{n+1}/H}_{M_{n}/H}=\mathrm{res}^{M/H}_{M_{n}/H}$, donc on d\'{e}duit un morphisme $\Phi: \mathrm{Gal}(M/H)\rightarrow \lim\limits_{\substack{\longleftarrow \\ n}}\mathrm{Gal}(M_{n}/H)$, qui est un isomorphisme car $(M_{n}/H)_{n \in \mathbb{N}}$ filtre $M/H$.
\fin
\vskip 2mm
\noindent
\tab Etant donnés un corps commutatif $\psi$-libre $k$ et un corps gauche $H$ contenant $k$ dans son centre, on voit que si pour tout système complet de groupes finis $(G_n,s_n)$, on sait trouver une tour $(L_n/k)$ réalisant ce système et telle que $H \otimes_{k} L_n$ reste un corps pour tout $n$, alors $H$ sera $\psi$-libre. C'est cette idée que nous allons utiliser pour montrer le théorème B. Nous exploitons à cet effet la construction présentée dans \cite{DD04} pour montrer la 
\vskip 2mm
\noindent 
{\bf Proposition 2.2.3.---} {\it Soit $k$ un corps commutatif contenant un corps $k_{0}$ qui est, soit r\'{e}el clos, soit hens\'{e}lien de caract\'{e}ristique r\'{e}siduelle nulle  et contenant $\mu_{\infty}$. Alors, pour tout syst\`{e}me complet de groupes finis $(G_{n},s_{n})_{n \in \mathbb{N}}$, il existe une tour d'extensions galoisiennes $(L_{n}/k(t))_{n \in \mathbb{N}}$ qui r\'{e}alise $(G_{n},s_{n})_{n \in \mathbb{N}}$ telle que $L_{n} \subset k((t))$ pour tout $n \in \mathbb{N}$.}
\vskip 2mm
\noindent
{\bf Preuve :} Il suffit de montrer le r\'{e}sultat lorsque $k=k_{0}$. En effet, supposons qu'il existe une tour $(F_{n}/k_{0}(t))_{n}$ qui r\'{e}alise $(G_{n},s_{n})_{n \in \mathbb{N}}$ telle que $F_{n} \subset k_{0}((t))$ pour tout $n \in \mathbb{N}$. Pour tout $n \in \mathbb{N}$, on pose  $L_{n}=F_{n}k$.  On obtient que la tour $(L_{n}/k(t))_{n \in \mathbb{N}}$ r\'{e}alise aussi $(G_{n},s_{n})_{n \in \mathbb{N}}$ et que l'on a $L_{n} \subset k((t))$ pour tout $n \in \mathbb{N}$. 
\vskip 2mm
\noindent
\tab Nous remarquons que, quitte \`{a} faire le changement de variables $t\mapsto 1/t$, il suffit de trouver une tour $(L_{n}/k(t))_{n \in \mathbb{N}}$ telle que $L_{n} \subset k((1/t))$ pour tout $n \in \mathbb{N}$. Pour la preuve, nous allons distinguer selon que le corps $k$ est r\'{e}el clos ou valu\'{e} hens\'{e}lien.
\vskip 2mm
\noindent 
{\bf -\underline{1er cas} : Le corps $k$ est r\'{e}el clos.}
\vskip 2mm
\noindent
Si $k=\mathbb{R}$, alors cela résulte de \cite[Remarque 2.2]{DD04}. Dans le cas g\'{e}n\'{e}ral, voir la preuve de \cite[Th\'{e}or\`{e}me 2.3]{DD04}.
\vskip 2mm
\noindent
{\bf -\underline{2\`{e}me cas} : Le corps $k$ est valu\'{e} hens\'{e}lien.} 
\vskip 2mm
\noindent
Notons $v$ la valuation sur $k$  et $(k_{v},v)$ le compl\'{e}t\'{e} de $(k,v)$. Par les constructions de \cite[Partie 3.1]{DD04}, il existe une suite croissante de parties finies non vides $(U_{n})_{n \in \mathbb{N}}$ de $k$ et une tour d'extensions galoisiennes $(\tilde{L}_{n}/k_{v}(t))_{n \in \mathbb{N}}$ de points de branchement\footnote{Etant donnée une extension finie galoisienne $E/K(t)$ $K$-régulière (c'est-à-dire $E \cap \overline{K}=K$), on dit que $t_0 \in \mathbb{P}^1(\overline{K})$ est un {\it point de branchement} de $E/K(t)$ si l'idéal $\langle t-t_0 \rangle$  est ramifié dans la clôture intégrale de $\overline{K}[t-t_0]$ dans $E\overline{K}$ (si $t_0=\infty$, $t-t_0$ doit être remplacé par $1/t$).} $U_{n}$ telles que $(\widetilde{L}_{n}/k_{v}(t))_{n \in \mathbb{N}}$ r\'{e}alise r\'{e}guli\`{e}rement (c'est-à-dire $\widetilde{L}_{n} \cap \overline{k_v}=k_v$) le syst\`{e}me $(G_{n},s_{n})_{n \in \mathbb{N}}$ et que pour tout $n \in \mathbb{N}$, la sp\'{e}cialisation\footnote{Etant donnés une extension finie galoisienne $E/K(t)$ $K$-régulière et $t_{0} \in \mathbb{P}^{1}(K)$, la {\it spécialisation} de $E/K(t)$ en $t_0$, notée $E_{t_0}/K$, est l'extension résiduelle en un idéal premier $\mathcal{P}$ au dessus de $\langle t-t_0 \rangle$. L'extension $E/K(T)$ étant galoisienne, le corps $E_{t_0}$ ne dépend pas du choix de $\mathcal{P}$.} $(\widetilde{L}_{n})_{\infty}$ de $\widetilde{L}_{n}/k_{v}(t)$ en $\infty$ est \'{e}gale \`{a} $k_{v}$. La preuve du \cite[Th\'{e}or\`{e}me 3.4]{DD04} montre que la tour $(\widetilde{L}_{n}/k_{v}(t))_{n \in \mathbb{N}}$ admet un mod\`{e}le  $(L_{n}/k(t))_{n \in \mathbb{N}}$ (c'est-à-dire $\widetilde{L}_{n}=L_{n}k_{v}$) qui r\'{e}alise r\'{e}guli\`{e}rement le syst\`{e}me complet $(G_{n},s_{n})_{n \in \mathbb{N}}$.  Pour tout $n \in \mathbb{N}$, puisque $\infty$ n'est pas un point de branchement de $\widetilde{L}_{n}/k_{v}(t)$, le point $\infty$ n'est pas non plus un point de branchement de $L_{n}/k(t)$. Vu que $(k,v)$ est hens\'{e}lien pour une valuation de rang $1$, on a $k_{v} \cap \overline{k}=k$. Or la sp\'{e}cialisation $(L_{n})_{\infty} $ de $L_{n}$ en $\infty$ est incluse à la fois dans $ \overline{k}$ et dans la sp\'{e}cialisation $(\widetilde{L}_{n})_{\infty}$ de $\widetilde{L}_{n}$ en $\infty$, donc $(L_{n})_{\infty}=k$.
\fin
\vskip 2mm
\noindent
\tab Nous avons maintenant les outils nécessaires pour établir la $\psi$-liberté d'un certain corps des fractions rationnels tordu : 
\vskip 2mm
\noindent
{\bf Théorème B'.---}{\it Soient $H$ un corps gauche de centre $k$ et $\sigma$ un automorphisme d'ordre fini de $H$. Supposons que $k^{\langle \sigma \rangle}$ contient un corps $k_{0}$ qui est, soit réel clos, soit hens\'{e}lien de caract\'{e}ristique r\'{e}siduelle nulle et contenant le groupe des racines de l'unité $\mu_{\infty}$. Alors, le corps $H(t,\sigma)$ est $\psi$-libre au sens de la définition 2.2.1.}
\vskip 2mm
\noindent
{\bf Preuve :}
Consid\'{e}rons un syst\`{e}me complet de groupes finis $(G_{n},s_{n})_{n \in \mathbb{N}}$. Par la proposition $2.2.3$, il existe une r\'{e}alisation $(L_{n}/k^{\langle \sigma \rangle}(t^{n}))_{n \in \mathbb{N}}$ du syt\`{e}me complet $(G_{n},s_{n})_{n \in \mathbb{N}}$ telle que $L_{n} \subset k^{\langle \sigma \rangle}((t^{n}))$ pour tout $n \in \mathbb{N}$. Par la proposition 2.1.2, l'algèbre $M_{n}=H(t,\sigma)\otimes_{k^{\langle \sigma \rangle}(t^{n})}L_{n}$ reste un corps pour tout $n \in \mathbb{N}$. La proposition 2.2.2 montre que $M/H=(\bigcup_{n}M_{n})/H$ est galoisienne extérieure, relativement algébrique à gauche et que la tour $(M_{n}/H(t,\sigma))_{n \in \mathbb{N}}$ est galoisienne finie, relativement algébrique à gauche réalisant $(G_{n},s_{n})_{n \in \mathbb{N}}$. Donc, $H(t,\sigma)$ est $\psi$-libre.
\fin
\vskip 2mm
\noindent
On voit que le théorème B' est équivalent au théorème B, comme remarqué précédemment.
\vskip 2mm
\noindent
{\bf Exemple d'application.---}Dans le théorème B, en prenant $H=\mathbb{H}$ le corps des quaternions de Hamilton (de centre $\mathbb{R}$) et $\sigma=\mathrm{Id}$, on obtient que le groupe $\widehat{F}_{\omega}$ est groupe de Galois sur $\mathbb{H}(t)$. De même, le groupe  $\widehat{F}_{\omega}$ est groupe de Galois sur $L(t,\sigma)$, où $L/\mathbb{C}((u))$ est une extension finie non triviale et $\sigma \in \mathrm{Gal}(L/\mathbb{C}((u)))$.
\bibliography{Biblio2}

\bibliographystyle{alpha}\vskip 4mm
\noindent
{\bf Angelot Behajaina}\\
{\sc Laboratoire de Math\'ematiques Nicolas Oresme, CNRS UMR 6139}\\
Universit\'e de Caen - Normandie\\
BP 5186, 14032 Caen Cedex - France\\
E-mail : angelot.behajaina@unicaen.fr  
\end{document}